\documentclass[12pt]{article}
\usepackage{amsmath,amssymb,amsfonts,amsthm,enumerate,extarrows}
\usepackage{eucal} 
\usepackage{xcolor}
\usepackage{mathabx}
\def\acts{\lefttorightarrow}

\newcommand{\cF}{\mathcal{F}}
\newcommand{\C}{\mathbb{C}}

\newcommand{\Z}{\mathbb{Z}}
\newcommand{\cO}{\mathcal{O}}
\newcommand{\cV}{\mathcal{V}}

\newcommand{\cM}{\mathcal{M}}
\newcommand{\cE}{\mathcal{E}}
\newcommand{\Pp}{\mathbb{P}^1}
\newcommand{\Pt}{\mathbb{P}^2}
\newcommand{\cH}{\mathcal{H}}
\newcommand{\cG}{\mathcal{G}}
\newcommand{\cB}{\mathcal{B}}
\newcommand{\cP}{\mathcal{P}}
\newcommand{\actson}{\acts}
\newcommand{\bmu}{{\boldsymbol \mu}}
\newcommand{\bnu}{{\boldsymbol \nu}}

\newcommand{\Gp}{\Gamma_+}

\newcommand{\ghat}{\hat{\mathfrak{g}}}
\newcommand{\bw}{\boldsymbol w}
\newcommand{\bv}{\boldsymbol v}
\newcommand{\ba}{\boldsymbol a}
\newcommand{\bb}{\boldsymbol b}
\newcommand{\bk}{\boldsymbol k}
\newcommand{\bl}{\boldsymbol l}
\newcommand{\bi}{\boldsymbol i}

\DeclareMathOperator{\diag}{diag}

\DeclareMathOperator{\End}{End}
\DeclareMathOperator{\Hilb}{Hilb}
\DeclareMathOperator{\Ext}{Ext}
\DeclareMathOperator{\rk}{rk}

\DeclareMathOperator{\ch}{ch}

\DeclareMathOperator{\Hom}{Hom}
\DeclareMathOperator{\Tr}{Tr}

\newtheorem{thm}{Theorem}

\newtheorem{lemma}{Lemma}
\newtheorem{cor}[lemma]{Corollary}
\newtheorem{prop}{Proposition}
\newtheorem{conj}{Conjecture}

\title{AGT and the Segal-Sugawara construction}
\author{Erik Carlsson}
\begin{document}
\maketitle
\begin{abstract}

The conjectures of Alday, Gaiotto and Tachikawa
\cite{AGT} and its generalizations 
have been mathematically formulated as
the existence of an action of
a $W$-algebra on the cohomology or $K$-theory of the 
instanton moduli space,
together with a Whitakker vector \cite{BFN,MO1,SchVas1}. 
However, the original conjectures
also predict intertwining properties with
the natural higher rank version
of the ``$\Ext^1$ operator'' which was previously studied
by Okounkov and the author 
in \cite{CO}, a result which is now sometimes referred 
to as AGT in rank one \cite{Alb,PSS}.
Physically, this corresponds to incorporating
matter in the Nekrasov partition functions,
an obviously important feature in the physical theory.
It is therefore of interest to 
study how the $\Ext^1$ operator relates to the 
aforementioned structures on cohomology in higher rank,
and if possible to find a formulation from which the 
AGT conjectures follow as a corollary.
In this paper, we carry out something analogous 
using a modified Segal-Sugawara construction
for the $\hat{\mathfrak{sl}}_2\C$ structure
that appears in Okounkov and Nekrasov's
proof of Nekrasov's conjecture \cite{NO1} for rank two.
This immediately implies the AGT identities
when the central charge is one, 
a case which is of particular interest for string theorists,
and because of the natural appearance of the Seiberg-Witten curve in this
setup, see for instance Dijkgraaf and Vafa \cite{DV}, as well as \cite{IKV}. 




\end{abstract}
\section{Introduction}

In \cite{AGT}, Alday, Gaiotto and Tachikawa
proposed a collection of identities of
explicit power series under a change of variables, 
each one associated a Riemann surface $\Sigma$
of genus $g=0,1$, and a list of $N$ marked points.
On one side of the equation are 
certain functions due to Nekrasov
from four-dimensional
supersymmetric gauge theory, which are defined as 
generating functions of 
equivariant localization integrals
over a smooth noncompact complex 
algebraic variety $\cM_{r,n}$ called
the instanton moduli space, with respect to a 
certain torus action defined below. 
Here $r,n$ are nonnegative integers called the
rank and instanton number, or the rank
and second Chern class in the description of $\cM_{r,n}$
as a moduli space of framed sheaves on $\C \Pt$.
See \cite{Nak1} for a thorough introduction
to this variety from several different descriptions.
On the other side
are correlation functions
from Liouville theory, defined in terms of
certain intertwiners of lowest weight
representations of the Virasoro algebra \cite{Tes}.
One of the variables in these identities is the central
charge $c$ of these representations on the Liouville side,
which is related to the torus parameters on the Nekrasov
side, as in \eqref{agtsubs}.

In an analogous way to the Hilbert-Chow map,
the space $\cM_{r,n}$ is a resolution 
of a singular space $\cM^0_{r,n}$ 
which contains the moduli space of
certain energy minimizing solutions to the 
four-dimensional Yang-Mills 
equation (instantons) modulo gauge symmetry, see
\cite{ADHM,DK,Nak1}. 
These solutions are important in gauge theory
because lowest energy solutions are expected to
dominate the contributions to certain path integrals.
In supersymmetric gauge theory, this approximation
is sometimes exact, through a principle called
supersymmetric localization, 
see \cite{Mir} for a general introduction.
Nekrasov's functions provide
a mathematically rigorous definition of such
integrals \cite{Nek1}.
A precise mathematical conjecture 
supporting this statement
due to Nekrasov relates certain limits of 
these quantities to the Seiberg-Witten free energy,
and was proved by Okounkov and Nekrasov in \cite{NO1}.
See Okounkov's ICM notes \cite{Ok1} for an exposition.

To physicists, the AGT relations reflect
the existence of a 
six dimensional superconformal 
quantum field theory on
\[X=\mathbb{R}^4 \times \Sigma,\]
predicted by the classification of superconformal
field theories. It is expected that one can then
recover the gauge theory partition function
by letting $\Sigma$ shrink away
through a process called compactification,
or dimensional reduction,
see \cite{Wit} as well as
\cite{Rod} for an introduction to the AGT conjectures.
Mathematically, one would like to extend the
AGT relations to more conceptual statements, 
which ideally would have implications to the physics.
One direction that has been carried out by several
authors \cite{BFN,MO1,SchVas1} in much greater generality 
for the case of pure (massless) gauge theory,
is to identify the cohomology of $\cM_{r,n}$ with a
Hilbert space, and construct the 
relevant conformal symmetry representations, 
in a similar manner to Nakajima's earlier
work \cite{Nak2}. In the case of
AGT in rank two, these structures include
an action of the Virasoro algebra on cohomology.

A more general version of the conjectures
would have to incorporate the presence of matter,
as is done in the original paper of AGT.
Mathematically, this may be encoded in the
$\Ext^1$ operator $W$ studied earlier by Okounkov and 
the author in the case of the Hilbert scheme of points
on a general surface \cite{CO}.
The main result of that paper identifies $W$ with
a ``vertex operator'' in terms of Nakajima's famous Heisenberg
operators, which completely calculates the
Nekrasov functions in rank one.
This result was used in \cite{Alb}, to
search for a basis of a Hilbert space 
on the Liouville side, which realizes the
$\Ext^1$ operator in the fixed point basis in rank two.
An obviously desirable and more direct approach
would be to formulate the cohomological structures
of the last paragraph in a way that produces
the intertwining properties with the $\Ext^1$
operator predicted by AGT.
Ideally, the original conjectures would then
follow as a corollary.



In theorem \ref{mainthm}, we
do something similar in the
setup of Nekrasov and Okounkov's proof
of Nekrasov's conjecture \cite{NO1}.
In that paper, the authors used a representation of
the affine special linear group 
$\mathfrak{\hat{sl}}_r \C$ to study the
``dual partition function,'' 
which is a sort of generating function of 
the original partition functions.
The main idea of this paper is to apply the Segal-Sugawara
construction to the
$\mathfrak{\hat{sl}}_r \C$ action for $r=2$,
and then to locate an extension of it that
intertwines the transformed $\Ext^1$ operator in the desired
way. The correct choice of action turns out not to be the most 
obvious candidate, as seen in proposition \ref{mainprop} below. 
This is a fortunate occurence because the usual
Segal-Sugawara construction would lead to 
non-irreducible representations
of the Virasoro algebra, and would derail the argument.

In corollary \ref{maincor},
we show that this implies
the AGT relations in the special case
where the central charge is set to one,
and $\Sigma$ is a genus one Riemann surface
with an arbitrary number of points removed.
Specializing the central charge
is simpler than the general ``refined'' case,
but is of particular interest to string theorists.
See for instance \cite{DV}, which predicts a proof of AGT
using matrix models for this case, as well as \cite{IKV}. 
Along the lines of \cite{DV},
the case of $c=1$ and $\Sigma=\C\Pp$ with four points
removed was studied by Morozov, Morozov and Shakirov 
\cite{MMS} using the Dotsenko-Fateev matrix model
and Selberg integrals as the starting point on
the Liouville side. The context of the dual-partition
function is also particularly
interesting because of the natural appearance of
the Seiberg-Witten curve.





\subsection{Acknowledgements}
The author acknowledges
support from the Simons Center for Geometry and Physics, 
Stony Brook University, as well as the
International Center for Theoretical Physics
at which some of the research for this paper was performed.
The author would also like to thank Andrey Smirnov for
providing some MAPLE code for numerically testing AGT on 
$\Pp$ with four points removed. This was
very helpful in sorting out the cases presented
in this paper.

\section{The AGT relations}
\subsection{Nekrasov functions}

We now recall some mathematical background to
define the Nekrasov functions. The definition
of the moduli space of framed torsion 
free sheaves on $\cM_{r,n}$ 
can be found in Nakajima's book \cite{Nak1} 
as well as \cite{HL}. For an introduction
to Nekrasov's functions and their mathematical 
and physical meaning, see 
\cite{Nek1,Ok1}.



Let $\left\{(x_0:x_1:x_2)\right\}$ coordinatize
the complex projective plane $\Pt=\C \Pt$,
and let  $\Pp_\infty \subset \Pt$ 
be the line at infinity defined by $x_0=0$. 
The moduli space of framed torsion free sheaves
is given set-theoretically by
\begin{equation}
\label{Mrn}
\cM_{r,n}=
\left\{(F,\Phi): \rk(F)=r, c_2(F)=n\right\}
\end{equation}
where $F$ is a torsion free sheaf on $\Pt$
which is locally free in a neighborhood of
$\Pp_\infty$, and $\Phi$ is a choice of isomorphism
\[\Phi : F\big|_{\Pp_\infty} 
\xlongrightarrow{\sim}
\cO^r_{\Pp_\infty}\]
called a framing at infinity.
This is a smooth noncompact complex algebraic 
variety of complex dimension $2rn$, also known as
the instanton moduli space.
When $r=1$, this space is isomorphic to the Hilbert
scheme of points in the complex plane, which 
as a set parametrizes certain ideals
in the ring $R=\C[x,y]$,
\begin{equation}
\label{M2H}
\cM_{1,n}\cong \Hilb_n \C^2 = 
\left\{I \subset R : \dim_{\C} R/I=n\right\}.
\end{equation}
The map is determined by restricting the sheaf 
to $\C^2= \Pt-\Pp_\infty$,
and using the framing $\Phi$ to obtain the inclusion map
to $R$.

There is a standard torus action on this space
\[G=T^2 \times T^r \lefttorightarrow \cM_{r,n},\]
defined as follows:
the action of $T^2$ is the induced one from
\begin{equation}
\label{gentorus}
T^2 = \C^* \times \C^* \lefttorightarrow \Pt,\quad
(z_1,z_2) \cdot (x_0:x_1:x_2)=(x_0:z_1^{-1} x_1 : z_2^{-1} x_2),
\end{equation}
by pullback of sheaves.
The action of $T^r$ is by
rotating the framing, i.e.
\[\bw \cdot \Phi=
\diag(\bw)\cdot \Phi,\quad
\bw=(w_1,...,w_r)\in T^r,\]
which commutes with the $T^2$ action.

The fixed points of this action are well known to be isolated, and
to correspond to $r$-tuples of partitions
$\boldsymbol \mu=(\mu^{(1)},...,\mu^{(r)})$ with
\[|\boldsymbol \mu|=\sum_{i} |\mu^{(i)}|=n.\]
Under the isomorphism \eqref{M2H} for $r=1$,
they correspond to the ideals
$I_\mu$ generated by all monomials $x^iy^j$ such that $(i,j)$
is not a box in the Young diagram of $\mu$, i.e.
$j\geq\mu_i$. In higher rank $r>1$,
the sheaf and framing $(\cF_{\bmu},\varphi_{\bmu})$
associated to such a diagram 
may be determined by the restriction to the plane
\[\cF_{\bmu} \big|_{\C^2} \cong F_{\bmu} :=
I_{\mu^{(1)}} \oplus \cdots \oplus I_{\mu^{(r)}}
\hookrightarrow R^r\]
together with the inclusion map on the right.

There is a bundle 
$\cE$ of rank $r(n_1+n_2)$
on $\cM_{r,n_1}\times \cM_{r,n_2}$ 
whose fiber over a point described by a pair of
sheaves $(\cF,\cG)$ is given by
\begin{equation}
\label{extbun}
\cE\big|_{\cF,\cG}=\Ext^{1}_{\Pt}(\cF,\cG(-\Pp_\infty)),
\end{equation}
The restriction of this bundle to the diagonal
when $n_1=n_2$ is well-known to be the tangent bundle.
In \cite{CO}, the $K$-theory class of the rank one case of this bundle
was defined for the Hilbert scheme of points on a general smooth
quasi-projective surface.
Its Euler class was proved to define an explicit vertex operator in
terms of Nakajima's Heisenberg operators, and was
further generalized to $K$-theory in \cite{CNO1}.
We have an action
\begin{equation}
\label{Gact}
G=T^2 \times T^r \times T^r \actson \cM_{r,n_0} \times \cM_{r,n_1}
\end{equation}
where $T^2$ acts diagonally, and the second and third tori
act on the framings on the first and second components.
This action lifts naturally to $\cE$ using the 
description \eqref{extbun}, making $\cE$ an equivariant
bundle.


Now define the torus characters of the fibers of this bundle
\[E_{\bmu,\bnu}(z_1,z_2,\bw,\bv):=\ch \cE_{\cF_\bmu,\cF_\bnu}
\in \Z[z_i^{\pm 1},w_i^{\pm 1},v_i^{\pm 1}]\]
where $\bw,\bv$ are elements of the first and second
$r$-dimensional torus in \eqref{Gact} respectively.
The answer can be expressed in terms of the answer for $r=1$ by
\[E_{\bmu,\bnu}(z_1,z_2,\bw,\bv)=
\sum_{i,j} w_i^{-1}v_j E_{\mu^{(i)},\nu^{(j)}}(z_1,z_2)\]
The rank one case may be calculated by restricting the
sheaves to the open subset $\C^2 \subset \Pt$,
\[E_{\mu,\nu}(z_1,z_2) = \chi_{R}(R,R)-\chi_{R}(I_\mu,I_\nu),\]
\begin{equation}
\label{chi}
\chi_{R}(F,G)= 
\sum_{i=0}^2 (-1)^i \ch \Ext^i_{R}(F,G) \in \Z((z_1,z_2))
\end{equation}
Now by the additivity of the equivariant Euler characteristic
on exact sequences, equation \eqref{chi} may be calculated in
any resolution of the modules $F,G$.
We then deduce the explicit formula
\begin{equation}
\label{chiexplicit}
\chi_R(I_\mu,I_\nu)=
\overline{M \ch I_\mu} \ch I_\nu
\end{equation}
where 
\[M=(1-z_1)(1-z_2),\quad \ch I_\mu= M^{-1}-Q_\mu,\quad
Q_\mu=\sum_{(i,j)\in \mu} z_1^i z_2^j\]
so that $\ch I_\mu$
is simply the torus character of $I_\mu$ as a vector space,
which lives in $\Z((z_1,z_2))$. The conjugation of 
one of these power series is determined by simply 
replacing $z_i=z_i^{-1}$ in its expression
as a rational function in $z_1,z_2$.

For instance, if $\mu,\nu=[1,1],[2,1]$, 
we would find that
\[\ch I_\mu=M^{-1}-1-z_2,\quad
\ch I_\nu = M^{-1}-1-z_1-z_2,\]
\[\chi_R(I_\mu,I_\nu)=
\left(1-(1-z_1^{-1})(1-z_2^{-1})(1+z_2^{-1})\right)\times\]
\[\left((1-z_1)^{-1}(1-z_2)^{-1}-1-z_1-z_2\right),\]
\begin{equation}
\label{Eex}
E_{\mu,\nu}(z_1,z_2)=
1+z_1^{-1}z_2+z_1z_2^{-2}+z_1^{-1}+z_2^{-1}.
\end{equation}
With some work, we can find
the explicit combinatorial expression
\begin{equation}
\label{hooks}
E_{\mu,\nu}(z_1,z_2)=\sum_{s\in\mu} z_1^{-a_{\mu}(s)-1}z_2^{l_\nu(s)}+
\sum_{s \in \nu} z_1^{a_\nu(s)}z_2^{-l_\mu(s)-1},
\end{equation}
where $a_\mu(s),l_\mu(s)$ are the arm and 
leg lengths of the box $s$
in the Young diagram of $\mu$ respectively. 
In the above expression,
these lengths may take negative values
if $s$ is not inside $\mu$. 
We may easily verify that this 
expression matches with example \eqref{Eex}.
See \cite{CO} for details.

We now recall the localization formula.
Suppose some variables $x_i$ are indentified as
$x_i=\exp(s_i)$, thinking of the $x_i$ as elements of 
a complex torus $\C^*$, and $s_i$ as elements of its Lie
algebra. We define
\begin{equation}
\label{euler}
e\left(\chi \right)=
\prod_i \left(\sum_j a_{ij} s_i\right)^{k_i},\quad
\chi=\sum_i k_i \prod_j x_j^{a_{ij}}
\end{equation}
which is the same as
the equivariant Euler characteristic of a representation with 
character $\chi$, viewed as an equivariant bundle over a point.
For the rest of the paper, we identify the following sets of variables:
\[z_\alpha=e^{t_\alpha},\quad w_\alpha=e^{a_\alpha},
\quad v_\alpha=e^{b_\alpha},\quad u_\alpha=e^{m_\alpha}\]
for any subscript $\alpha$, which might in fact be a pair of indices,
i.e. $w_\alpha=w_{ij}$ as we have below.
For instance, we would have
\[e\left(z_1u_1^{-1}-2w_{12}u_2\right)= \frac{t_1-m_1}{(a_{12}+m_2)^2}.\]
Let us also set
\[e_m\left(\chi\right)=e(e^m\chi)=
\prod_i \left(m+\sum_j a_{ij} t_i \right)^{k_i}.\]

The equivariant localization formula for a smooth, 
projective
torus equivariant variety $T^d \actson X$ 
with isolated fixed points,
and a cohomology class $\gamma \in H_T(X)$, states that
\begin{equation}
\label{loc}
\int_X \gamma = \sum_{p\in X^T} \frac{i_p^* \gamma}{e(T_p X)}
\end{equation}
where $i_p$ is the inclusion of the fixed point $p$, and
the integral denotes the proper pushforward map to a point,
see \cite{AB,SS}.
If $z_i=e^{t_i}\in T^d$ are the torus variables, then this expression
is written as an element of $\C(t_1,...,t_d)$, but in fact it must
reside in $\C[t_1,...,t_d]$, the equivariant cohomology of a point, 
implying some cancelation. 
If $X$ is not compact, then \eqref{loc} may be taken as a definition,
extending integration to a functional satisfying
\[\int_X \gamma =\int_Y \pi_*(\gamma),\quad \gamma:X\rightarrow Y\]
for proper maps $\pi$. For some cohomology classes $\gamma$
in a suitable completion of the cohomology ring $H_T(X)$,
this definition
can be shown to coincide with the usual integration of differential
forms vanishing rapidly at infinity even for 
some noncompact manifolds $X$.

Now fix a positive integer $N$, and
for any letter $x$, let $\tilde{x}=(x_1,...,x_N)$
denote an $N$-tuple of indexed variables. For a bold letter
denoting an $r$-tuple symbols such as $\boldsymbol a$ above, let
\[\tilde{\boldsymbol a}=
\left({\boldsymbol a}_1,...,{\boldsymbol a}_{N}\right),\quad
\boldsymbol a_i= (a_{i1},...,a_{ir})\]
For partitions we will use the superscript, since the subscript is reserved:
\[\tilde{\boldsymbol \mu}=
\left\{\boldsymbol \mu^{(1)},...,\boldsymbol \mu^{(N)}\right\},\quad
\bmu^{(i)}=(\mu^{(i1)},...,\mu^{(ir)})\]

The Nekrasov functions are defined by
\[Z (t_1,t_2,\tilde{\boldsymbol a},\tilde{m},\tilde{q})=
\sum_{n_1,...,n_N} q^{n_1+\cdots+n_N}
\int_{\cM_{r,n_1}\times \cdots \times \cM_{r,n_N}}
\prod_{i=1}^N \pi_{i,i+1}^* (\cE)=\]
\begin{equation}
\label{nekdef}
\sum_{\tilde{\boldsymbol \mu}}
\prod_{i=1}^N q_i^{|\bmu^{(i)}|}
\frac{w_{\bmu^{(i)},\bmu^{(i+1)}}(t_1,t_2,\boldsymbol a_i,\ba_{i+1},m_i)} 
{w_{\bmu^{(i)},\bmu^{(i)}}(t_1,t_2,\ba_i,\ba_i,0)},
\end{equation}
where $\pi_{i,j}$ is the projection onto the $i$th and $j$th factor,
\[w_{\bmu,\bnu}(t_1,t_2,\ba_i,\ba_j,m)=
e_m\left(E_{\bmu,\bnu}(z_1,z_2,
{\boldsymbol w}_i,{\boldsymbol w}_j)\right),\]
and we identify $N+1$ with $1$. The sum is over all $N$-tuples
of $r$-tuples of partitions, not just those of a fixed norm.
For $N=1$, this is the ``instanton part'' of the partition functions originally
written down by Nekrasov in \cite{Nek1} 
in the form of contour integrals. For a mathematical introduction to the meaning
of these integrals, we refer the reader to 
Okounkov's ICM notes \cite{Ok1}.

For example, for $r=2,N=1$ we would have
\[Z(t_1,t_2,(a,-a),m,q)=\]
\[1+\left(\frac{(m-2a)(m-t_2)(m-t_1)(m+2a-t_1-t_2)}{2at_1t_2(t_1+t_2-2a)}-
\right.\]
\begin{equation}
\label{Zex}
\left.\frac{(m+2a)(m-t_2)(m-t_1)(m-2a-t_1-t_2)}{2at_1t_2(t_1+t_2+2a)}\right)q+O(q^2)
\end{equation}
up to first order in $q$.

\subsection{Conformal Blocks}
The second side of the AGT relations are formed by the conformal
blocks of the Virasoro algebra. 
See \cite{MMMM1} 
for many useful calculations in the context of AGT.

Recall the commutation relations for the Virasoro algebra
given by
\begin{equation}
\label{vir}
[L_{m},L_{n}] = (m-n)L_{m+n}+\delta_{m,n}\frac{m^3-m}{12}K,
\end{equation}
where $K$ is central.
Let $M_{h}=M_{h,c}$ denote the Verma module of the Virasoro 
algebra of level $h$ and central charge $c$,
so that the Cartan subalgebra generated by $L_0,K$
acts on its lowest weight (vaccuum) vector $v_{\emptyset}$ by
\[X v_{\emptyset}=\Lambda_{h,c}(X)v_\emptyset,\quad
\Lambda_{h,c}(L_0)=h,\quad
\Lambda_{h,c}(K)=c,\]
whereas positive generators $L_k$ for $k>0$ annihilate it.
There is a basis of $M_{h}$ indexed by partitions given by
\[v_{\mu}=L_{-\mu} v_\emptyset,\quad
L_{-\mu}= L_{-\mu_1}\cdots L_{-\mu_l}.\]
We will also write $v_{\mu,h}$ if we want to stress
that $v_{\mu,h}$ is an element of $M_{h}$.

There an inner product on $M_h$ called
the Shapovalov form defined by
\[(v_\emptyset,v_\emptyset)_h=1,\quad
(L_{-n} v_{\mu},v_{\nu})_h=(v_{\mu},L_{n} v_\nu)_h.\]
We denote its matrix elements by
\[K_{\mu,\nu}(h)=(v_{\mu},v_{\nu})_h\]
For instance, the entries $K_{\mu,\nu}(h)$
for $\mu,\nu\in \{[2],[1,1]\}$ are given by
\[\left(\begin{array}{cc}
4h+c/2 & 6h \\
6h & 8h^2+4h
\end{array}\right).\]

There is an intertwining operator called the 
\emph{Liouville vertex operator}
\[\cV^{h}_{k_2,k_1}(x) : M_{k_1} \rightarrow M_{k_2}\]
satisfying
\begin{equation}
\label{LCO1}
\left(\cV^{h}_{k_2,k_1}(x) v_{\emptyset,k_1},v_{\emptyset,k_2}\right)_{k_2}=x^{k_2-h-k_1},
\end{equation}
\begin{equation}
\label{LCO}
[L_n,\cV^{h}_{k_2,k_1}(x)]=
\left(h(n+1)x^n+x^{n+1}\partial_x\right)\cV^{h}_{k_2,k_1}(x).
\end{equation}
We will find it more convenient to use the normalization
\[V^{h}_{k_2,k_1}(x)=x^{k_1+h-k_2} \cV^{h}_{k_2,k_1}(x),\]
so that the vaccuum matrix element is one. 
Strictly speaking, 
this operator is only defined as a \emph{field}, i.e.
a formal power series 
\[\cV^{h}_{k_2,k_1}(x) \in \Hom(M_{k_1},M_{k_2})[[x^{\pm 1}]],\]
with the property that each matrix element is a formal 
Laurent series, see \cite{FBZ}.
However, 
the full Liouville vertex
operator may be defined to act on a larger Hilbert space of 
functions on a noncompact space,
see \cite{Tes} for an introduction.

We set
\[S_{\mu,\nu}(k_1,h,k_2)=S_{k_1,h,k_2}(v_{\mu},v_\nu)=
\left(\cV^{h}_{k_2,k_1}(x) v_{\mu,k_1},
v_{\nu,k_2}\right)_{k_2}\Big|_{x=1}.\]
%
For instance, we have
\[S_{[1],[1]}(k_1,h,k_2)=
S_{k_1,h,k_2}(L_{-1}v_{\emptyset},L_{-1}v_{\emptyset})=\]
\[S_{k_1,h,k_2}(v_{\emptyset},L_1 L_{-1} v_{\emptyset})+
(k_1+h-k_2-1)S(v_{\emptyset},L_{-1}v_\emptyset)=\]
\begin{equation}
\label{Sexample}
2k_2+(k_1+h-k_2-1)(k_2+h-k_1).
\end{equation}
More complicated terms will of course involve the central charge $c$.

Let $d=L_0-h$ so that $d v_\mu = |\mu| v_{\mu}$.
Let
\[\cB(c,\tilde{h},\tilde{k},\tilde{q})=
\Tr q^d V^{h_{N}}_{k_1,k_{N}}(x_N)\cdots
V^{h_{2}}_{k_3,k_{2}}(x_{2})V^{h_1}_{k_2,k_1}(x_1)=\]
\[\sum_{\tilde{\mu},\tilde{\nu}}
q_1^{|\mu^{(1)}|}\cdots q_N^{|\mu^{(N)}|}
\frac{S_{\nu_1,\mu_2}(k_1,h_1,k_2)\cdots S_{\nu_N,\mu_1}(k_{N},h_{N},k_1)}
{K_{\mu_1,\nu_1}(k_1) \cdots
K_{\mu_N,\nu_N}(k_N) }\]
where 
\begin{equation}
\label{x2q}
q=x_1=q_1\cdots q_N,\quad
x_ix_{i+1}^{-1}=q_{i+1}
\end{equation}

\subsection{The AGT conjecture}


The AGT relations for a genus 1 Riemann
surface with $N$ punctures state that
\begin{conj}(AGT \cite{AGT})
Let $\ba_i=(a_i,-a_i)$ for $1\leq i\leq N$. Then we have that
\begin{equation}
\label{agteq}
Z(t_1,t_2,\tilde{\ba},\tilde{m},\tilde{q})=
Z'(t_1,t_2,\tilde{m},\tilde{q})\cB(c,\tilde{k},\tilde{h},\tilde{q})
\end{equation}
under the substitution
\[c=1+6\frac{(t_1+t_2)^2}{t_1t_2},\quad
k_i=\frac{(t_1+t_2)^2-4a_i^2}{4t_1t_2},\]
\begin{equation}
\label{agtsubs}
h_i=\Delta_{m_i,m_i},\quad \Delta_{m,n}=\frac{m(t_1+t_2-n)}{t_1t_2},
\end{equation}
and
\[Z'(t_1,t_2,\tilde{m},\tilde{q})=
(q;q)_\infty^{2\Delta_{m_1,m_1}+\cdots+2\Delta_{m_N,m_N}-1}\times\]
\[\prod_{i<j} (x_ix_j^{-1};q)^{2\Delta_{m_i,m_j}}_\infty 
(x_i^{-1}x_j q;q)^{2\Delta_{m_j,m_i}}_\infty\]
where the $x_i$ are determined by \eqref{x2q}, and
$(x;q)_\infty=\prod_{i \geq 0} (1-xq^i).$
\end{conj}

For instance, for $N=1$ we would have
\[Z'(t_1,t_2,m,q)=(q;q)_\infty^{2\Delta_{m,m}-1}=
1+\left(1-\frac{2m(t_1+t_2-m)}{t_1t_2}\right)q+\cdots\]
and
\[\cB(c,k,h)=1+S_{[1],[1]}(k,h,k)K^{-1}_{[1],[1]}(k)q+\cdots=\]
\[1+\frac{h^2-h+2k}{2k}q+\cdots\]
using our example \eqref{Sexample} at $k_1,k_2=k$.
Using \eqref{Zex}, we can check that
\[Z(t_1,t_2,(a,-a),m,q)=Z'(t_1,t_2,m,q)\cB(c,k,h,q)\]
to first order in $q$, after the substitution
\[k=\frac{(t_1+t_2)^2-4a}{t_1t_2},\quad
h=\frac{m(t_1+t_2-m)}{t_1t_2}.\]
The central charge $c$ does not appear until higher order.

\section{Vertex operators}
\subsection{The infinite wedge representation}


Let $\Lambda=\Lambda^{\infty/2}$ denote the infinite wedge
representation, which we will now briefly summarize 
and refer to \cite{BO,Kac} for details.
It has a basis labeled by partitions which we
write as
\[v_{\mu}=v_{\mu_1}\wedge v_{\mu_2-1}\wedge 
v_{\mu_3-2}\wedge \cdots\]
The possible sequences
\[(i_1,i_2,i_3,...)=(\mu_1,\mu_2-1,\mu_3-2,...)\] 
that can appear are precisely the strictly decreasing
sequences of integers such that the number of
entries with $i_k>0$ equals the number of $i_k\leq 0$.
This space can be 
thought of as the $\infty/2$ exterior
power of the vector space $\C\cdot \Z$ with basis 
$v_i$ indexed by the
integers. The wedge product of sums
of basis vectors $v_i$ should be distributed and
sorted with signs in the usual way. For instance,
\[v_{2}\wedge (v_{1}+v_{4})\wedge v_{-2}\wedge \cdots=\]
\[(v_{2}\wedge v_{1}\wedge v_{-2}\wedge \cdots) -
(v_{4}\wedge v_{2}\wedge v_{-2}\wedge \cdots)=
v_{[2,2]}-v_{[4,3]}.\]

There is an action of the Lie algebra of 
the infinite-dimensional
general linear group $\mathfrak{gl}(\C\cdot \Z)$ defined by
\[\rho'(E_{ij}) \cdot v_{i_1}\wedge v_{i_2}\wedge \cdots=\]
\[\sum_{k} v_{i_1} \wedge \cdots \wedge v_{i_{k-1}} \wedge
E_{ij} v_{i_k} \wedge v_{i_{k+1}} \wedge \cdots\]
where
\[E_{ij} v_k=\begin{cases}
v_i & \mbox{$j=k$} \\
0 & \mbox{otherwise}
\end{cases}\]
This extends to a projective representation
defined by
\begin{equation}
\label{infwedge}
\rho(E_{ij}) = \begin{cases}
\rho'(E_{ij}) & \mbox{ if $i\neq j$ or $i=j>0$} \\
\rho'(E_{ij})-Id & \mbox{ if $i=j\leq 0$}
\end{cases}
\end{equation}
which has the commutation relations
\[[\rho(E_{ij}),\rho(E_{kl})]=\rho([E_{ij},E_{kl}])+
\epsilon_{ijkl} Id,\]
\begin{equation}
\label{EijCR}
\epsilon_{ijkl}=\begin{cases}
1 & \mbox{$i=l \leq 0$ and $j=k\geq1$}\\
-1 & \mbox{$i=l \geq 1$ and $j-k \geq 0$}\\
0 & \mbox{otherwise}
\end{cases}
\end{equation}

The projective representation extends to some infinite
sums of the elementary matrices $E_{ij}$ which are useful
for defining the action of Ka\c{c}-Moody algebras
on $\Lambda$. For instance, the element
\[d = \sum_{i\in \Z} i \rho\left(E_{ii}\right)\]
becomes a finite sum when applied to any vector,
and is determined by
\[d\cdot v_{\mu}=|\mu|v_{\mu}.\]
We also define the action of the infinite-dimensional
Heisenberg Lie algebra
\[\alpha_n=\sum_{i\in \Z} \rho\left(E_{i,i+n}\right)\]
satisfying
\begin{equation}
\label{alphacr}
[\alpha_m,\alpha_n]=m\delta_{m,n}Id.
\end{equation}

There is an isomorphism from the polynomial algebra
\[\C[\alpha_{-1},\alpha_{-2},...] \cong
\Lambda,\]
given by simply applying the polynomial on the left
to the vacuum $v_{\emptyset}$, so the image of the monomials would be
\[\alpha_{\lambda} = \alpha_{-\lambda_1}\cdots \alpha_{-\lambda_l} \cdot v_\emptyset\]
for a Young diagram $\lambda$ of length $l$. To determine the inverse map
amounts to finding coefficients in the expansion
\[v_\mu=\sum_{\lambda} c_{\lambda,\mu} \alpha_{\lambda}.\]
It turns out that the $c_{\lambda,\mu}$ are precisesly the 
coefficients of the expansion of the Schur polynomial $s_\mu$
in the power sum basis $p_\lambda$ 
\cite{Mac}.


Now for each $m$, we have the following 
well-known \emph{vertex operator},
\begin{equation}
\label{VO}
\Gamma^{(m)}(x) = \Gamma_-^m(x)\Gamma_+^{-m}(x^{-1}) \in 
\End\left(\Lambda\right)[[x^{\pm 1}]],
\end{equation}
where
\[\Gamma_{\pm}^m(x)= \exp\left(
m\sum_{k>0} \frac{x^{k}\alpha_{\pm k}}{k}\right).\]
We will also use the odd and even parts, given by
\[\Gamma^{(m)}_{s}(x)=\Gamma^{m}_{s,-}(x)\Gamma^{-m}_{s,+}(x^{-1})\]
where $s=e,o$, and
\[\Gamma_{e,\pm}(x)=
\exp\left(m\sum_{k>0,even} \frac{x^{k}\alpha_{\pm k}}{k}\right),\quad
\Gamma_{o,\pm}(x)=
\exp\left(m\sum_{k>0,odd} \frac{x^{k}\alpha_{\pm k}}{k}\right)\]
so that $\Gamma^{(m)}(x)=
\Gamma^{(m)}_e(x) \Gamma^{(m)}_o(x)$.
We will also write
\begin{equation}
\label{eodecomp}
\gamma\Gamma^{(m)}(x)\gamma^{-1}=
\Gamma^{(m)}_e(x) \otimes \Gamma^{(m)}_o(x)
\end{equation}
where 
\[\gamma : \Lambda \rightarrow
\C[\alpha_{-1},\alpha_{-2},...] \cong\]
\[\C[\alpha_{-2},\alpha_{-4},...] \otimes
\C[\alpha_{-1},\alpha_{-3},...] 
=:
\Lambda_e \otimes \Lambda_o.\]

By exponentiating the commutation relations 
\eqref{alphacr}, we can determine that
\begin{equation}
\label{Omega}
\Gamma_{s,+}^{m}(x)\Gamma_{s,-}^{n}(y)=
\Omega_s(x,y)^{mn}\Gamma_{s,-}^{n}(y)\Gamma_{s,+}^{m}(x),
\end{equation}
where
\[\Omega_e(x,y)=\frac{1}{\sqrt{1-x^2y^2}},\quad
\Omega_o(x,y)=\frac{\sqrt{1-x^2y^2}}{(1-xy)}\]
so that
\[\Omega_e(x,y)\Omega_o(x,y)=\Omega(x,y)=\frac{1}{1-xy}.\]
We also have the easier relation
\begin{equation}
\label{gcrq}
q^d \Gamma_{s,\pm}(x)= \Gamma_{s,\pm}(xq^{\mp 1})q^d.
\end{equation}
\subsection{The principal vertex operator construction}

Now let $\ghat$ 
denote the affine Lie algebra $\mathfrak{\hat{sl}}_2 \C$, 
whose underlying vector space is given by
\[\ghat=\mathfrak{\hat{sl}}_2\C=
\mathfrak{sl}_2(\C[t,t^{-1}])+\C d' + \C K.\]
Let
\[e=\left(\begin{array}{cc}0 & 1 \\ 0 & 0\end{array}\right),\quad
h=\left(\begin{array}{cc}1 & 0 \\ 0 & -1\end{array}\right),\quad
f=\left(\begin{array}{cc}0 & 0 \\ 1 & 0\end{array}\right)\]
denote the standard generators of $\mathfrak{sl}_2\C$. For
each $i\in \Z$ and $a\in \mathfrak{sl}_2\C$, let 
\[a(x)=\sum_{k\in \Z} a_{-k}x^k \in \ghat[[x^{\pm 1}]],\quad
a_k=a\cdot t^{-k}\in \ghat\]
The commutators are given by
\[[a_m,b_n]=[a,b]_{m+n}+m\delta_{m,-n}K,\]
\begin{equation}
\label{sl2hat}
[d',a_k]=ka_{k-1},\quad [K,\ghat]=0.
\end{equation}

There is an action of $\ghat$ on $\Lambda$
which is induced from the action
\[\mathfrak{sl}_2\left(\C[t,t^{-1}]\right) \acts \C^2 [t,t^{-1}]
\cong \C\cdot \Z\]
where the isomorphism is defined by
\[  (t^k,0) \mapsto v_{2k},\quad
(0,t^k) \mapsto v_{2k-1}.\]
Explicitly, we have
\[e_i \mapsto \sum_k E_{2k,-1+2i+2k},\quad
h_i \mapsto \sum_k E_{2k,2i+2k}-E_{-1+2k,-1+2i+2k},\]
\begin{equation}
\label{sl2wedge}
f_i \mapsto \sum_k E_{-1+2k,2i+2k},\quad
2d'+\frac{1}{2}h_0 \mapsto d,\quad K\mapsto 1,
\end{equation}
where it is understood the the elementary 
matrices $E_{ij}$ act via 
the projective representation $\rho$.
This representation is not irreducible, but the span of
$\ghat$ applied to $v_\emptyset$ is isomorphic to the basic
representation $\Lambda_0$.


The principal vertex operator construction gives an
additional description of the action \eqref{sl2wedge},
which explicitly identifies it as the space
$\Lambda_o \subset \Lambda$ defined above.
It is described by
\[2d' \mapsto d-A_0/2,\quad
K \mapsto 1,\]
\begin{equation}
\label{pvoc}
h_i\mapsto A_{2i},\quad
2e_i\mapsto \alpha_{2i-1}-A_{2i-1},\quad
2f_i\mapsto \alpha_{2i+1}+A_{2i+1}
\end{equation}
where
\[2A(x)=2\sum_i A_{-i}x^i = \Gamma^{(2)}_o(x)-1.\]


\subsection{The Segal-Sugawara construction}

We now recall some facts about the Segal-Sugawara construction,
for which we refer to \cite{FBZ}.
Let $V$ be a representation
of an affine Lie algebra $\ghat$ with central charge $c$.
The Segal-Sugawara construction produces operators
$L_n$ in terms of the generators of $\ghat$ such that
\begin{enumerate}
\item The vector space $V$ becomes a representation
of the Virasoro algebra with central charge
\[c'=\frac{c \dim \mathfrak{g}}{c+h^{\vee}},\]
where $h^\vee$ is the dual Coxeter number.
In particular, if $\mathfrak{g}=\mathfrak{sl}_2$ and
$V=\Lambda_0$, we would have 
$c=1,h^\vee=2,\dim \mathfrak{g}=3$, so that $c'=1$.
\item The Virasoro algebra interwines the action of 
$\ghat$ in the desired manner, coming from the action
of automorphism of the circle on 
$Map(S^1,\mathfrak{g})$ by
precomposition, i.e.
\[[L_m,a(x)]= x^{1+m}\partial_x a(x)
\quad a \in \mathfrak{g}.\]
Furthermore, $L_0$ coincides with the differential $d'$.
\end{enumerate}

In the case $\mathfrak{g}=\mathfrak{sl}_2\C$, the 
construction has the form
\begin{equation}
\label{sugsl2}
L_k=\frac{1}{12}\sum_{i\in \Z} 
:2e_if_{k-i}+2f_ie_{k-i}+h_ih_{k-i}:
\end{equation}
where the ``normal ordering'' symbol means
\[:a_ib_j:=
\begin{cases}
a_ib_j & \mbox{if } i\leq 0 \\
b_ja_i & \mbox{otherwise.}
\end{cases}\]
Combining this formula with \eqref{sl2wedge},
we arrive at an action of the Virasoro algebra
with central charge $1$ on $\Lambda$
which preserves $\Lambda_o$.
In fact, there is a 
family of actions of $\cV ir$ parametrized
by a number $s$ given by
\begin{equation}
\label{Ls}
L_{k,s}=L_k+sh_k+s^2\delta_{k,0}.
\end{equation}
It is straightforward to verify that these also 
satisfy \eqref{vir}. 
For integer values of $s$ these are
the translates of the original action by the
translation subgroup of the affine 
Weyl group of $\hat{\mathfrak{sl}}_2$.

We have a decomposition of $\Lambda_o$
as a representation over $L_{k,s}$ for
any $s$ as follows.
By the Ka\c{c} character formula, we have that
\begin{equation}
\label{kacsl2}
\Tr_{\Lambda_o} y^{h_0}q^{d'}=
\sum_{k\in \Z} y^{2k}q^{k^2}(q;q)_{\infty}^{-1},\quad
(x;q)_\infty=\prod_{i\geq0} (1-xq^i).
\end{equation}
Since the Virasoro action commutes with $h_0$,
we find that $\Lambda_o$ decomposes as
\[\Lambda_o \cong \bigoplus_k V_k\]
where $V_k$ is the eigenspace of $h_0$ with eigenvalue $2k$,
and contains a unique up to scalar lowest
eigenvector $v_{k}$ of $d'$ with eigenvalue $k^2$, given by
$v_k=v_{\mu}$ where
\begin{equation}
\label{vk}
\mu=\begin{cases}
[2k,2k-1,...,1] & k\geq 0 \\
[-2k-1,-2k-2,...,1] & k<0
\end{cases}
\end{equation}
We find that
\[L_{0,s} v_{k}= (d'+sh_0+s^2)v_k=(k+s)^2v_k,\]
giving rise to a map
\begin{equation}
\label{V2M}
M_{(k+s)^2} \rightarrow V_k.
\end{equation}
Furthermore, if $s$ is in the range where
$M_{(k+s)^2}$ is irreducible, then this map
is injective. Using \eqref{kacsl2}, we find
that 
\[\Tr_{M_{(k+s)^2}} q^{L_{0,s}}=
\Tr_{V_k} q^{L_{0,s}}=
q^{(k+s)^2}(q;q)_\infty^{-1}\]
so the map is also an isomorphism in that range.


The following proposition 
will be used to prove our main theorem:

\begin{prop}
\label{mainprop}
We have that
\begin{equation}
\label{indcom}
[L_{k,1/4},\Gamma^{(2m)}_o(x^{1/2})]=
\left(m^2kx^k+x^{k+1}\partial_x\right)\Gamma^{(2m)}_o(x^{1/2})
\end{equation}
%
\end{prop}

\begin{proof}


We begin with the case $m=1$.
In this case, using \eqref{pvoc}, we have that
\[\Gamma^{(2)}_o(x^{1/2})=
-2x^{1/2}e(x)+2h(x)+2x^{-1/2}f(x)+1=\]
\[x^{h_0/4} a(x) x^{-h_0/4}+1,\quad
a=-2f+2h+2f.\]
Then we have
\[[L_{k,1/4},\Gamma^{(2)}_o(x^{1/2})]=
[L_{k}+h_k/4+1/16,x^{h_0/4} a(x) x^{-h_0/4}+1]=\]
\[x^{h_0/4}\left(kx^k+x^{k+1}\partial_x\right)a(x)x^{-h_0/4}+\]
\[1/4\cdot \left(x^kx^{h_0/4}[h,a](x)x^{-h_0/4}+4kx^k\right)=\]
\[kx^kx^{h_0/4}a(x)x^{-h_0/4}+
x^{k+1}\partial_x \left(x^{h_0/4}a(x)x^{-h_0/4}\right)+kx^k=\]
\[\left(kx^k+x^{k+1}\partial_x\right)\Gamma^{(2)}_o(x^{1/2}).\]

We now suppose now that \eqref{indcom} holds for some $m,n$,
and prove that it holds for $m+n$.
Using \eqref{Omega}, we find that
\[\Gamma^{(2(m+n))}_o(x^{1/2})=
\lim_{y\rightarrow x}
\Gamma_{-,o}^{2m}(x^{1/2})
\Gamma_{-,o}^{2n}(y^{1/2})
\Gamma_{+,o}^{-2m}(x^{-1/2})
\Gamma_{+,o}^{-2n}(y^{-1/2})=\]
\[\lim_{y\rightarrow x}
A(x,y)
\Gamma^{(2m)}_o(x^{1/2})
\Gamma^{(2n)}_o(y^{1/2}),\quad
A(x,y)=\Omega_o(x^{-1/2},y^{1/2})^{4mn}\]
We then have
\[[L_{k,1/4},\Gamma^{(2m+2n)}_o(x^{1/2})]=\]
\[\lim_{y\rightarrow x}
A(x,y)
(km^2x^k+kn^2y^k+x^{k+1}\partial_x+y^{k+1}\partial_y)\]
\[\Gamma^{(2m)}_o(x^{1/2})\Gamma^{(2n)}_o(y^{1/2})=\]
\[\lim_{y\rightarrow x}
A(x,y)
(km^2x^k+kn^2y^k+x^{k+1}\partial_x+y^{k+1}\partial_y)\]
\[A(x,y)^{-1}
\Gamma_{-,o}^{2m}(x^{1/2})
\Gamma_{-,o}^{2n}(y^{1/2})
\Gamma_{+,o}^{-2m}(x^{-1/2})
\Gamma_{+,o}^{-2n}(y^{-1/2})=\]
\[(k(m^2+n^2)x^k+x^{k+1}\partial_x)
\Gamma^{(2m+2n)}_o(x^{1/2})+\]
\[\left(\lim_{y\rightarrow x}A(x,y)
(x^{k+1}\partial_x+y^{k+1}\partial_y)
A(x,y)^{-1}\right)
\Gamma^{(2m+2n)}_o(x^{1/2}).\]
The limit in parentheses in the last line 
can be seen to equal $2kmnx^k$, which combines with the first expression
to establish \eqref{indcom}.


\end{proof}

\section{Application to AGT}

\subsection{Nekrasov functions for $t_1+t_2=0$}
We now explain how the infinite wedge
representation was applied to Nekrasov functions
in \cite{NO1}.

First, we restrict our torus action by specializing
\begin{equation}
\label{cytorus}
z_1=z,\quad z_2=z^{-1},\quad t_1=t,\quad t_2=-t,\quad z=e^t
\end{equation}
Now substitute \eqref{cytorus} into \eqref{chiexplicit} to get
\[E_{\bmu,\bnu}(z,{\boldsymbol w},{\boldsymbol v})=
E_{\bmu,\bnu}(z,z^{-1},{\boldsymbol w},{\boldsymbol v})=
\sum_{i,j=1}^r w_i^{-1}v_j E_{\mu,\nu}(z),\]
\begin{equation}
\label{cEz}
E_{\mu,\nu}(z) =
\chi_{\emptyset,\emptyset}(z)-\chi_{\mu,\nu}(z),\quad
\chi_{\mu,\nu}(z)=f_{\mu}(z^{-1})f_{\nu}(z)
\end{equation}
where $f_\mu(z)\in\C(z)$ is the rational function whose Laurent
series about infinity is given by
\[f_\mu(z)=\sum_{i\geq 1} z^{\mu_i-i+1}=
\sum_{i=1}^{\ell(\mu)} z^{\mu_i-i+1}+\frac{z^{-\ell(\mu)-1}}{(1-z^{-1})}.\]
Notice that the expression in \eqref{cEz} must be a polynomial,
implying some cancellation.

For the rest of the paper we will assume without any 
additional loss of information that $t=1$
and write
\[Z(\tilde{\ba},\tilde{m},\tilde{q})=
Z(1,-1,\tilde{\ba},\tilde{m},\tilde{q}),\]
\[Z'(\tilde{m},\tilde{q})=Z'(1,-1,\tilde{m},\tilde{q}),\quad
\cB(\tilde{h},\tilde{k})=\cB(1,\tilde{h},\tilde{k}).\]
Set
\begin{equation}
\label{wdef}
w_{\bmu,\bnu}(\ba,\bb,m)=e_m\left(E_{\bmu,\bnu}(z)\right),\quad
w_{\mu,\nu}(m)=e_m\left(E_{\mu,\nu}(z)\right)
\end{equation}
By the symmetry $z \leftrightarrow z^{-1}$, there
is a polynomial $w_{\bmu}({\boldsymbol a})$ such that
\begin{equation}
\label{hbmu}
w_{\bmu}(\ba) :=
w_{\bmu,\bmu}({\boldsymbol a},{\boldsymbol a},0)=
(-1)^{r|\bmu|} w_{\bmu}({\boldsymbol a})^2
\end{equation}
If $r=1$, the $w_{\mu}$ is simply the product of the hook lengths of $\mu$.

For instance, we would have
\[f_{[2,1]}(z)= z^{2}+1+\frac{z^{-2}}{1-z^{-1}},\quad
f_{[1,1]}(z)= z+1+\frac{z^{-2}}{1-z^{-1}}.\]
This would give
\[E_{[1,1],[2,1]}(z)= \frac{1}{(1-z)(1-z^{-1})}-f_{[1,1]}(z^{-1})f_{[2,1]}(z)=\]
\[z^3+z+1+z^{-1}+z^{-2},\]
which agrees with \eqref{Eex} under the specialization
\eqref{cytorus}.
We then have
\begin{equation}
\label{emex}
w_{[1,1],[2,1]}(m)=(m + 3) (m + 1) m (m - 1) (m - 2).
\end{equation}

Now let $\cH^r$ denote the complex vector space with basis
vectors given by $u_{\bmu}$, where $\bmu$ is an $r$-tuple of partitions,
and let $\cH=\cH^1$. We define an operator
\[W_{\ba,\bb}^{(m)}(x) \in \End(\cH^r)[[x^{\pm 1}]],\]
\begin{equation}
\label{Wdef}
\left(W_{\bb,\ba}^{(m)}(x) u_{\bmu},u_{\bnu}\right)=
(-1)^{r|\bmu|}x^{|\bnu|-|\bmu|}
\frac{w_{\bmu,\bnu}(\ba,\bb,m)}{w_{\bmu}(\ba)w_{\bnu}(\bb)}.
\end{equation}
If $r=1$ we will simply write $W^{(m)}(x)$.
We then have
\begin{equation}
\label{tragt}
Z(\tilde{\ba},\tilde{m},\tilde{q})=
\Tr_{\cH^r} q^d 
W^{(m_N)}_{\ba_1,\ba_N}(x_N)\cdots
W^{(m_{2})}_{\ba_3,\ba_{2}}(x_{2})
W^{(m_1)}_{\ba_2,\ba_1}(x_1)
\end{equation}
where the $x_i$ are related to the $q_i$ by \eqref{x2q}.

Consider the following isomorphism:
\[\iota : \cH \rightarrow \Lambda,\quad
u_\mu \mapsto v_\mu\]
The following result was proved
in \cite{NO1}, and was extended to
the general action \eqref{gentorus} 
(and in fact, to a general
smooth quasiprojective surface)
by Okounkov and the author in \cite{CO}. 
A further extension of this theorem to $K$-theory
may be found in \cite{CNO1}. 
The author has also described
a very short proof for the 
more general $K$-theoretic version,
but for the specialized 
action \eqref{cytorus} in \cite{C1}.
\begin{prop}
\label{voprop}
We have that
\begin{equation}
\label{VOeq}
\iota W^{(m)}(x) \iota^{-1}=
\Gamma^{(m)}(x).
\end{equation}
\end{prop}

For instance, dropping the power of $x$, we
find that
\[(\Gamma^{(m)} v_{[1,1]},v_{[2,1]})=
(\Gp^{-m}v_{[1,1]},\Gp^{m} v_{[2,1]}).\]
Next, we get
\[\Gp^{-m}v_{[1,1]}=
\Gp^{-m} v_{1}\wedge
\Gp^{-m} v_{0}\wedge
\Gp^{-m}v_{-2}\wedge \cdots=\]
\[\left(v_{1}-mv_{0}+\frac{m(m-1)}{2}v_{-1}+\cdots \right)
\wedge \]
\[\left(v_{0}-mv_{-1}+\cdots\right)
\wedge \left(v_{-2}+\cdots\right)\wedge \cdots=\]
\[v_{[1,1]}-mv_{[1]}+\frac{m^2+m}{2}v_{\emptyset},\]
and
\[\Gp^{m} v_{[2,1]}=
\Gp^m v_{2}\wedge
\Gp^m v_{0}\wedge
\Gp^mv_{-2}\wedge \cdots=\]
\[\left(v_{2}+mv_{1}+\frac{m(m+1)}{2}v_{0}+
\frac{m(m+1)(m+2)}{6}v_{-1}+\cdots\right)\wedge\]
\[\left(v_{0}+mv_{-1}+
\cdots \right)\wedge
\left(v_{-2}+\cdots\right)\wedge \cdots=\]
\[v_{[2,1]}+mv_{[2]}+mv_{[1,1]}+m^2v_{[1]}+
\frac{m^3-m}{3}v_{\emptyset}.\]
Taking the inner product of the two yields
\[\frac{m(m-1)(m-2)(m+3)(m+1)}{6}\]
which agrees with \eqref{emex}, \eqref{Wdef},
and \eqref{VOeq}.


%
We now explain how to apply proposition \ref{voprop}
to higher rank, which was used by Okounkov and 
Nekrasov to compute
the \emph{dual partition function} in \cite{NO1}.
Consider the function
\[\cB_r : \cP^r \times \Z^r \rightarrow \cP \times \Z\]
which associates to an $r$-tuple a \emph{blended partition}
\[(\boldsymbol\mu,\boldsymbol k)=
(\mu^{(1)},...,\mu^{(r)};k_1,...,k_r)\mapsto
(\mu,k),\]
where $\mu,k$ are the uniquely determined by
the property that
\begin{equation}
\label{blend}
\left\{\mu_{i}-i+1+k\right\}_{i\geq 1}=
\bigcup_{j=1}^r \left\{r
\left(\mu^{(j)}_i-i+1+k_j\right)-j+1\right\}_{i\geq 1}
\end{equation}
and $k=|{\boldsymbol k}| =k_1+\cdots+k_r$.
The norms are related by
\begin{equation}
\label{blendnorm}
|\mu|=r|\bmu|+d_{\boldsymbol k},\quad
d_{\boldsymbol k}=\frac{r-1}{2}\sum_{i} k_i^2+\frac{r+1-2i}{2}k_i-
\sum_{i<j} k_ik_j.
\end{equation}
It is straightforward to see that this map is 
bijective. If we have $|{\boldsymbol k}|=0$ then we obtain 
another isomorphism
\[\beta_{\boldsymbol k} : \cH^r \rightarrow \cH,\quad
u_{\bmu} \mapsto u_{\mu},\quad \cB_r(\bmu,{\boldsymbol k})=(\mu,0)\]




\begin{prop}
\label{dpfprop}
We have
\begin{equation}
\label{dpfeq}
\beta_{\boldsymbol l}^{-1}
W^{(rm)}(x)
\beta_{\boldsymbol k}
= c_{{\boldsymbol k},{\boldsymbol l},m}
x^{d_{\boldsymbol l}-d_{\boldsymbol k}}
W_{\bl-\bi/r,\bk-\bi/r}^{(m)}(x^r)
\end{equation}
where $\bi=(1,2,...,r)$, and $c_{{\boldsymbol k},{\boldsymbol l},m}$ is a constant.
\end{prop}

\begin{proof}
If
\[(\mu,0)=\cB_r(\bmu,{\boldsymbol k}),
\quad (\nu,0)=\cB_r(\boldsymbol\mu,\boldsymbol l),\]
then we have
\[f_{\mu}(z)=\sum_{j=1}^r z^{k_j-j+1}f_{\mu^{(j)}}(z^r),\]
and therefore
\[E_{\mu,\nu}(z)-E_{\mu^0,\nu^0}(z)=
E_{\boldsymbol\mu,\boldsymbol\nu}(z^r,z^{r\bk-\bi},z^{r\bl-\bi})\]
where
\[(\mu^0,0)=
\cB_r({\boldsymbol \emptyset},{\boldsymbol k}),\quad
(\nu^0,0)=\cB_r({\boldsymbol\emptyset},\boldsymbol l)\]
are the blended $r$-tuples of empty partitions.

It follows that
\[w_{\mu,\nu}(rm)\sim
e_{rm}\left(E_{\boldsymbol\mu,\boldsymbol\nu}(z^r,z^{r\bk-\bi},z^{r\bl-\bi})\right)=\]
\[r^{2r(|\bmu|+|\bnu|)}w_{\bmu,\bnu}({\boldsymbol k}-{\boldsymbol i}/r,
{\boldsymbol l}-{\boldsymbol i}/r,m)\]
where $\sim$ means the left side is a multiple of the
right side by a constant that does not depend on $\bmu,\bnu$.
This also implies that
\[w_{\mu}(rm) \sim r^{2r|\bmu|}w_{\bmu}({\boldsymbol k}-{\boldsymbol i}/r,
{\boldsymbol l}-{\boldsymbol i}/r,m)\]
Inserting both of these and \eqref{blendnorm} 
into \eqref{Wdef} yields the result up to the sign, which
is straightforward to determine.
\end{proof}

\subsection{The main theorem}

We may now state our main theorem.

\begin{thm}
\label{mainthm}
Let $\bk=(k,-k),\bl=(l,-l)$, and consider the composition
\[\gamma_{\bk} : \cH^2 
\xlongrightarrow{\beta_{\bk}} \cH
\xlongrightarrow{\iota} \Lambda
\xlongrightarrow{\gamma} \Lambda_e \otimes \Lambda_o.\]
Then
\begin{enumerate}[a)]
\item The map $\gamma_{\bk}$ is injective, and its image is
\label{firstthmpart}
$\Lambda_e\otimes V_k$. We may therefore write
\[\gamma_{\bk} : \cH^2 \rightarrow \Lambda_e \otimes V_k\]
and refer to its inverse $\gamma_{\bk}^{-1}$.
\item We have that
\label{secondthmpart}
\begin{equation}
\label{mainthmeq}
\gamma_{\bl} W^{(m)}_{\bl-\bi/2,\bk-\bi/2}(x)
\gamma_{\bk}^{-1} =
\Gamma^{(2m)}_{e}(x^{1/2}) \otimes V^{m^2}_{(l+1/4)^2,(k+1/4)^2}(x),
\end{equation}
where $V^{h}_{k_2,k_1}(x)$ is
the Liouville vertex operator defined above
using the map \eqref{V2M} with $s=1/4$,
which is an isomorphism for this value.
\end{enumerate}
\end{thm}

\begin{proof}
First, it is easy to see that the image of $\gamma_{\bk}$ is
the eigenspace of $h_0$ with eigenvalue $2k$, from which
part \ref{firstthmpart} follows.


By propositions \ref{voprop} and \ref{dpfprop}, we have
\begin{equation}
\label{gammaVO}
\gamma_{\bl} W^{(m)}_{\bl-\bi/2,\bk-\bi/2}(x)
\gamma_{\bk}^{-1}=
\Gamma_{e}^{(2m)}(x^{1/2}) \otimes
x^{d_{\bk}/2-d_{\bl}/2}\Gamma_{o}^{(2m)}(x^{1/2})
\end{equation}
Now we have
\[d_{\bl}/2-d_{\bk}/2=(l+1/4)^2-(k+1/4)^2,\]
so we must show that
\[x^{-m^2}\Gamma_{o}^{(2m)}(x^{1/2})\]
satisfies \eqref{LCO}. But 
this follows from proposition \ref{mainprop}.
To determine that the vaccuum expectation is one
on both sides of \eqref{gammaVO}, 
it is enough 
to check values on the lowest weight vectors,
and notice that the partition 
in \eqref{vk} is precisely the blended partition
\[(\mu,0)=\cB(\emptyset,\emptyset;k,-k).\]

\end{proof}

We now have
\begin{cor}
\label{maincor}
The AGT relations hold for $t_1+t_2=0$.
\end{cor}
\begin{proof}
For an $N$-tuple of integers $\tilde{k}$, let
\[\bk_i=(k_i,-k_i),\quad a_i=k_i+1/4,\quad \ba_i=(a_i,-a_i)\]
so that
\[W^{(m)}_{\ba_i,\ba_j}(x)=
W^{(m)}_{\bk_i-\bi/2,\bk_j-\bi/2}(x),\]
by subtracting the constant 3/4 from 
both the entries of $\ba_i,\ba_j$.
It suffices to prove the claim for these values because
the coefficients of $q_i$ are rational functions,
which are determined by their values at these points.

Now apply part \ref{secondthmpart} of the theorem
to 
\eqref{tragt} to get
\[Z(\tilde{\ba},\tilde{m},\tilde{q})=
\Tr q^d W^{(m_N)}_{\bk_1-\bi/2,\bk_N-\bi/2}(x_N)\cdots
W^{(m_1)}_{\bk_{2}-\bi/2,\bk_{1}-\bi/2}(x_1)=\]
\[\left(\Tr_{\Lambda_e} q^{d/2} \Gamma^{(2m_N)}_e(x_N^{1/2})
\cdots \Gamma^{(2m_1)}_e(x_1^{1/2})\right)\times\]
\[\left(\Tr q^d V^{m_N^2}_{(k_1+1/4)^2,(k_N+1/4)^2}(x_N)\cdots
V^{m_1^2}_{(k_2+1/4)^2,(k_1+1/4)^2}(x_1)\right).\]
The second factor is by definition
$\cB(\tilde{k},\tilde{h},\tilde{q})$ under the
change of variables in \eqref{agteq},
so it remains to show that the first factor
equals $Z'(\tilde{m},\tilde{q})$.
This can be calculated using the
commutation relations \eqref{Omega} and \eqref{gcrq}.
We will verify it for $N=1$, leaving the general
case as an exercise:
\[\Tr_{\Lambda_e} q^{d/2} \Gamma^{(2m)}_e(x^{1/2})=
\Tr q^{d/2} \Gamma_{-,e}^{2m}(x^{1/2})
\Gamma_{+,e}^{-2m}(x^{-1/2})=\]
\[\Tr \Gamma_{-,e}^{2m}((xq)^{1/2})q^{d/2}
\Gamma_{+,e}^{-2m}(x^{-1/2})=\]
\[\Tr q^{d/2} \Gamma_{+,e}^{-2m}(x^{-1/2})
\Gamma_{-,e}^{2m}((xq)^{1/2})=\]
\[(1-q)^{2m^2}\Tr q^{d/2} \Gamma_{-,e}^{2m}((xq)^{1/2})
\Gamma_{+,e}^{-2m}(x^{-1/2})=\]
\[\cdots\]
\[(q;q)_\infty^{2m^2} \Tr q^{d/2} \Gamma_{+,e}(x^{-1/2})\]
Since $\Gamma_+$ is unitriangular with respect to the 
degree grading, we get
\[(q;q)_\infty^{2m^2} \Tr_{\Lambda_e} q^{d/2}=(q;q)_\infty^{2m^2-1}=
Z'(m,q).\]
\end{proof}

\bibliography{sugbib}{}

\end{document}